\documentclass[10pt]{article}
\usepackage{amsthm,amsmath,amsfonts}
\usepackage{fullpage,enumerate}
\usepackage{tikz}
\usepackage[colorlinks,citecolor=blue,urlcolor=blue]{hyperref}
\usepackage{hypernat}
\usepackage[notref,notcite]{showkeys}
\usepackage{graphicx}

\newcommand{\la}{\lambda}

\begin{document}

\numberwithin{equation}{section}

\title{Bacterial persistence: a winning strategy?}

\bigskip
\author{  Olivier Garet, R\'egine Marchand, \\{\it Universit\'e de Lorraine, Nancy} \and Rinaldo B. Schinazi\\
{\it University of Colorado, Colorado Springs}}
\maketitle

\bigskip
{\bf Abstract.}  It has long been known that antibiotic treatment will not completely kill off a  bacteria population. For many species a small fraction of bacteria is not sensitive to antibiotics. These bacteria are said to persist. Recently it has been shown that persistence is not a permanent state and that in fact a bacterium can switch back and forth between persistent and non persistent states.  We introduce two stochastic models for bacteria persistence. In both models there are mass killings of non persistent bacteria at certain times. The first model has deterministic killing times and the second one has random killing times. Both models suggest that persistence may be a successful strategy for a wide range of parameter values. 
\bigskip

{\bf Key words:} bacteria persistence, stochastic model, branching process, random environment

\bigskip

{\bf AMS Classification:} 60K35

\bigskip

{\bf 1. Introduction}

Since at least Bigger (1944) it is known that antibiotic treatment will not completely kill off a  bacteria population. For many species a small fraction of bacteria is not sensitive to antibiotics.
These bacteria are said to be 'persistent'. It turns out that persistence is not a permanent state. In fact a bacterium can switch back and forth between a 'normal' (i.e. non persistent) state and a 'persistent' state.  In a normal state a bacterium can reproduce but is killed by an antibiotic attack.  In the persistent state a bacterium does not reproduce (or very seldom) but is resistant to antibiotics, see Balaban et al. (2004), Kussell et al. (2005) and Lewis (2007). Persistent bacteria may play an important role in treatment failure and in the appearance of bacteria with inherited resistance to antibiotics (a major health concern), see Levin (2004),  Wiuff et  al. (2005) and Levin and Rozen (2006).  

In this paper we consider two models for this phenomenon. 
In both models there are mass killings of non persistent bacteria at certain times.
The first model has deterministic killing times and the second one has random killing times. Both models are random in that giving birth to a new bacterium and switching back and forth between persistent and non persistent states all happen randomly. It is believed that switching is random for several species of bacteria, for the biology of persistence see Dubnau and Losick (2006) and Graumann (2006).

For the deterministic killing time model we assume that at certain fixed times all the normal bacteria are killed off and only the persistent ones survive. We will show that for all parameter values  the bacteria population has a positive probability of surviving provided the time interval between two successive killings is larger than a certain critical value. We will also show that this critical value is not very sensitive to even large variations in parameter values.  This seems to suggest that persistence is a successful strategy for a wide range of parameter values. We also show a similar behavior for the model where the killing occurs at random times. For this second model  we concentrate on the particular case for which time intervals between two killings is an i.i.d. sequence of exponentially distributed random variables.

\bigbreak
{\bf 2. Deterministic killing times}

Consider the following continuous time model. Bacteria can be in one of two states. We think of state 1 as being the normal (i.e. non persistent) state and state 2 as being the persistent state. Note that the normal state is vastly predominant in the population. 
A bacterium in state 1 is subject to two possible transitions. It can give birth at rate $\lambda$ to another state 1 individual or it can switch to state 2 at rate $a$. A bacterium in state 2 has only one possible transition. It can switch to state 1 at rate $b$. Moreover, the bacteria in state 1 can be killed in the following way.
Let $T_i=iT$ for $i\geq 1$ where $T$ is a fixed positive constant. This defines a sequence of killing times $T_1,T_2,\dots$. At each killing time $T_i$ all the bacteria in state 1 are killed  but the bacteria in state 2 are unaffected. 

The model mimics the observed behavior. State 1 bacteria multiply until disaster strikes and then they all die. State 2 bacteria cannot give birth but persist under disasters. Hence, state 2 bacteria ensure survival through disasters but cannot give birth. Our hypothesis that bacteria in state 2 cannot give birth or die is not unrealistic. See Table 1 in Kussel et al. (2005).

The main question we are concerned with in this paper is for which parameter values do the bacteria survive? The following result answers this question.
\medskip

{\bf Theorem 1. }{\sl For any $a>0$, $b>0$ and $\la>0$ there is a critical value $T_c$ such that the bacteria survive forever with positive probability if and only if $T>T_c$.}

\medskip
The exact value of $T_c$ depends on the parameters $a$, $b$ and $\la$. However, $T_c$ varies very little with $a$ and $b$, see Figure 1. The value of $T_c$ is more sensitive to variations of $\la$. See Figure 2.

\bigbreak

\includegraphics[width=6cm]{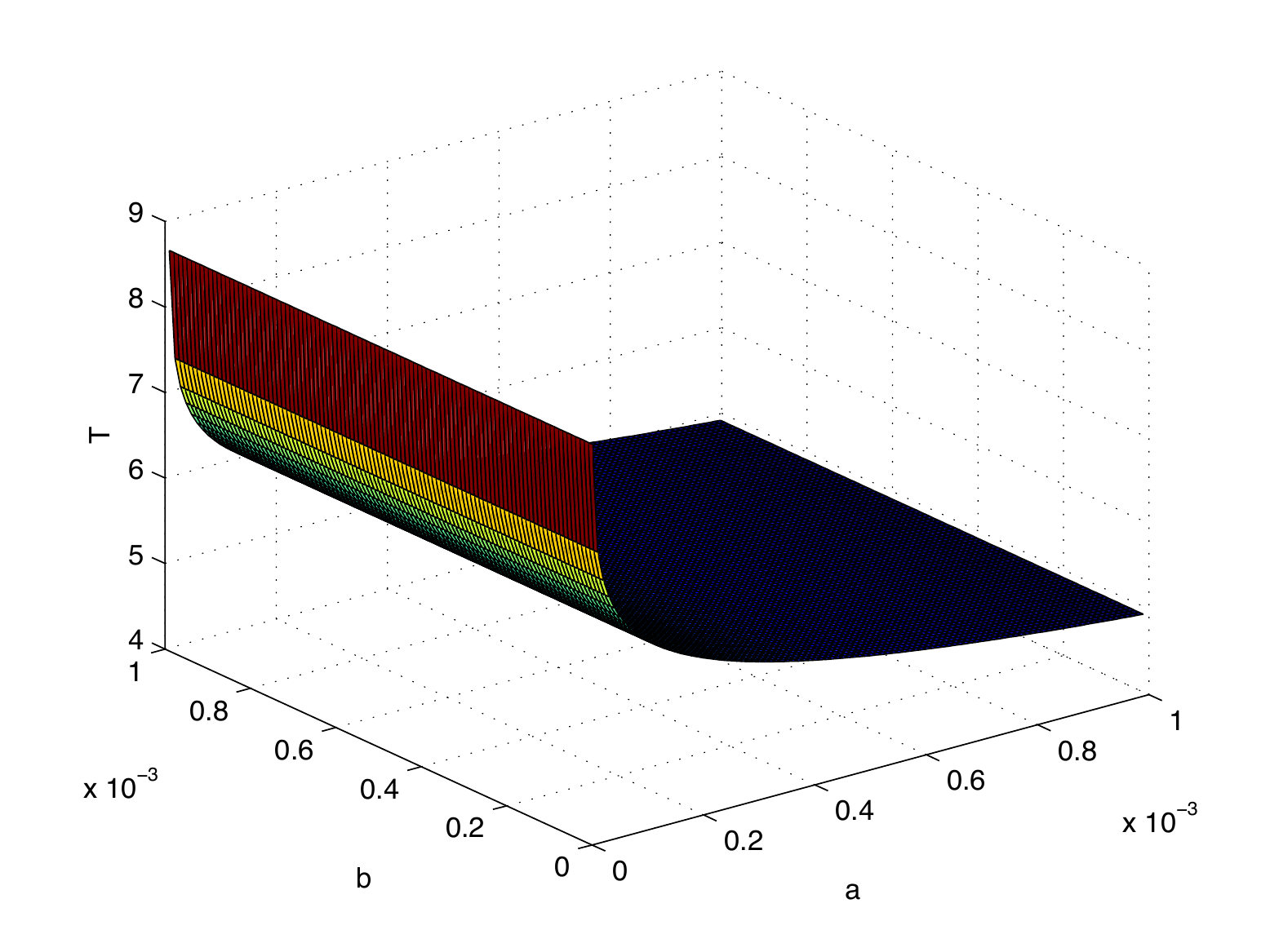}

{\bf Figure 1.} This is the critical value $T_c$ as a function of $a$ and $b$. The parameters $a$ and $b$ vary from $10^{-6}$ to $10^{-3}$ and $\la$ is fixed at 2. 

\bigbreak

\includegraphics[width=6cm]{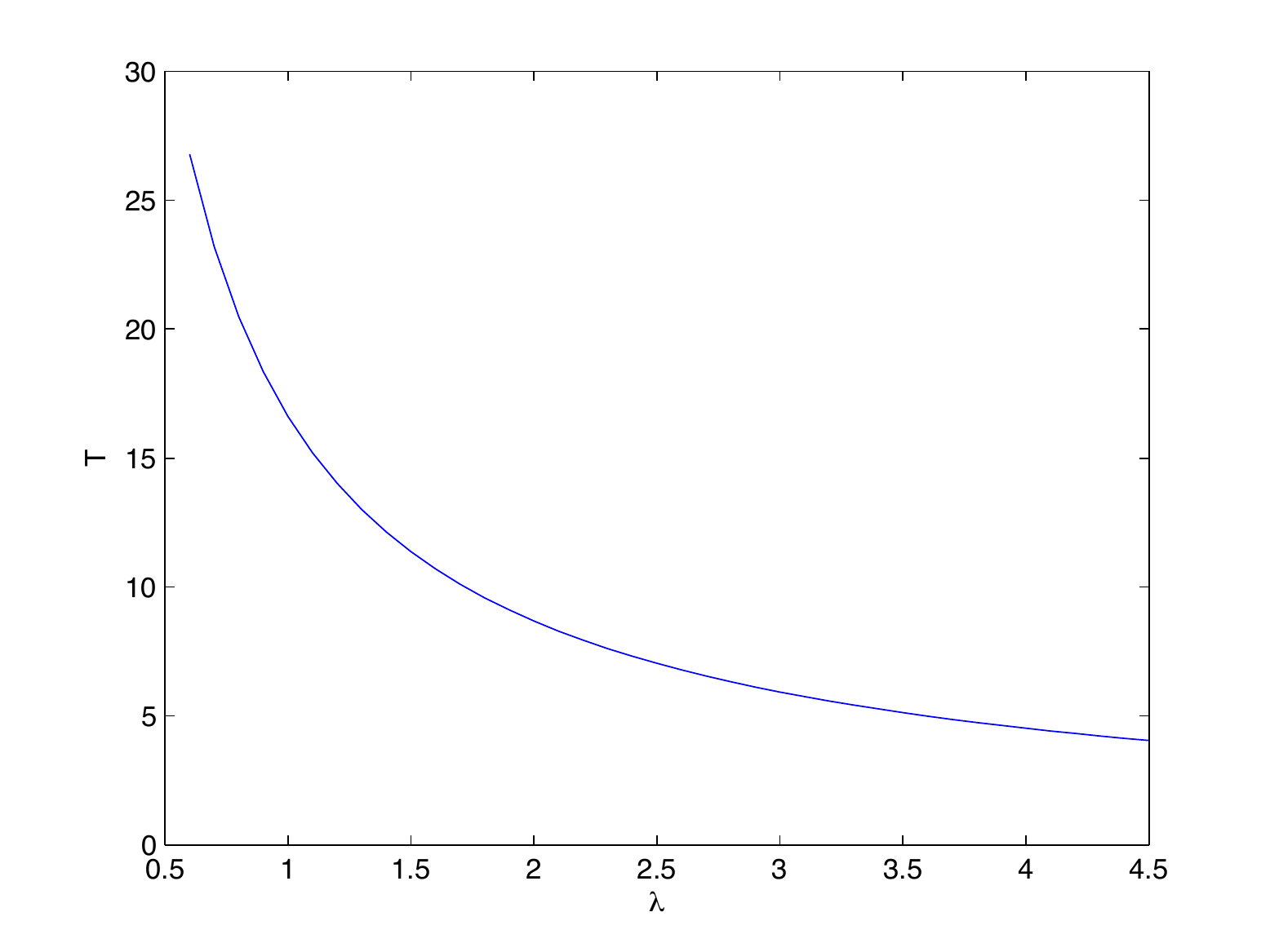}

{\bf Figure 2.} This is the critical value $T_c$ as a function of $\la$. We set $a=10^{-6}$ and $b=10^{-3}$. 
\bigbreak
{\bf 2. Random killing times}

We now turn to random killing times.
The growth and switching back and forth rates $\la$, $a$ and $b$ are the same as before. The only difference with the model above is that the killings now occur at random times.
We consider a sequence $(T_k)_{k\geq 1}$ of i.i.d. random times which is also independent of everything else in the model. At each $T_k$ we kill all the bacteria in state 1 and the state 2 bacteria are unaffected. Many choices are possible for the distribution of the random times. We choose to concentrate on the exponential distribution. That is, we assume that $(T_k)_{k\geq 1}$ is an i.i.d. sequence with probability density $\delta\exp(-\delta t)$ where $\delta>0$.
\bigbreak

{\bf Theorem 2. }{\sl  Let $a>0$, $b>0$ and $\la>0$. 
There exists $\delta_c>0$ such that:
\begin{itemize}
\item if $\delta<\delta_c$, then the model with random killing times has a positive probability of surviving,
\item if $\delta \ge \delta_c$, then the model with random killing times always dies out.
\end{itemize}
%
%
}
Note that the expected killing time is $1/\delta$. It plays the role of $T$ in the previous model. 
\bigbreak

Balaban et al. (2004) and Kussell et al. (2005) have used differential equations to model persistence. In their model they use a certain set of parameter values when the environment is a 'growth' environment and a different set of values when the environment is a 'stress' environment. They use the model to compute long term fitnesses for persistent and normal states. Our Theorem 1 shows that even for very small values of $a$, $b$ and $\la$ the bacteria may survive provided the interval between two killing times is long enough is consistent with their findings. However, their model unlike ours is very sensitive to parameter variations, see their Discussion. 

Garner et al. (2007) also use differential equations to model persistence but they are closer to our approach in that they kill all the normal bacteria at certain fixed times. They use their model to compute the ideal fraction of persistent bacteria in a population. Note that Gardner et al. (2007) also consider deterministic and random (exponentially distributed) killing times. \bigbreak

{\bf 3. Proof of Theorem 1}

We start our process with finitely many bacteria. 
We now define an auxiliary discrete time stochastic process $Z_n$, $n\geq 0$. We wait until the first killing time $T$ and we let $Z_0$ be the number of bacteria in state 2 at time $T$. If $Z_0=0$ we set $Z_i=0$ for $i\geq 1$. If $Z_0\geq 1$ then we wait until the second killing time $2T$ and let $Z_1$ be the number of  state 2 bacteria at time $Z_1$. More generally, for any $k\geq 1$ let $Z_k$ be the number of state 2 individuals at the $(k+1)$-th  killing time $(k+1)T$. 

We claim that the process $(Z_k)_{k\geq 0}$ is a Galton-Watson process. This can be seen using the following argument. Each state 2 bacterium present at time $T_1=T$ starts a patch of bacteria.  Note that in order for the patch to get started we first need the initial state 2 bacterium to switch to state 1.  At the second killing time $T_2$ all the individuals in state 1 are killed and we are left with $Z_1$ individuals in state 2. If $Z_1=0$ the bacteria have died out. If $Z_1\geq 1$ then each 2 present at time $T_1$ starts its own patch. Each patch will live between times $T_2$ and $T_3$. These patches are independent and identically distributed. At time $T_3$ each 2 (if any) starts a new patch and so on.  
This construction shows that $(Z_k)_{k\geq 0}$ is a Galton-Watson process. Moreover, the bacteria population survives forever if and only if $Z_k\geq 1$ for all $k\geq 0$. For if $Z_k=0$ for some $k\geq 1$ this means at the corresponding killing time $T_{k+1}$ there are no type 2 individuals and no type 1 individuals either.  That is, the bacteria have died out. If on the other hand $Z_k\geq 1$ for all $k\geq 0$ the bacteria will survive forever. Now, it is well known that a Galton-Watson process survives if and only if $E(Z_1|Z_0=1)>1$, see for instance Schinazi (1999). That is, the analysis can be done starting with a single state 2 bacterium. Hence, the problem of survival for the bacteria population is reduced to computing the expected value $E(Z_1|Z_0=1)$. We do this computation next.

For $t<T_1$ let $x(t)$ be the expected number of type 1 bacteria and $y(t)$ be the expected number of type 2 bacteria at time $t$, starting at time 0 with a single type 2 and no type 1. We have for $h>0$
\begin{eqnarray*}
x(t+h)-x(t) & = & \lambda h x(t)-a hx(t)+bhy(t)+o(h), \\
y(t+h)-y(t) & = & ahx(t)-bhy(t)+o(h).
\end{eqnarray*}
This yields the following system of differential equations
\begin{eqnarray*}
{d\over dt}x(t) &= &(\lambda-a) x(t)+by(t), \\
{d\over dt}y(t) &= & a x(t)-by(t).
\end{eqnarray*}
This is a linear system with constant coefficients. The corresponding matrix 
$$
A=\begin{pmatrix} 
\la-a& b\\
a& -b\\
\end{pmatrix}
$$
has two real distinct eigenvalues $\nu_1$ and $\nu_2$
$$\nu_1=\frac{-a - b + \la + \sqrt{\Delta}}{ 2}$$
$$\nu_2=\frac{-a - b + \la - \sqrt{\Delta}}{2},$$
where
$$\Delta=(a+b-\la)^2+4b\la>(a+b-\la)^2.$$
Note that the determinant of $A$ is $-\la b<0$. Hence, $\nu_1>0$ and $\nu_2<0$. 
A standard computation yields the solution of the system of differential equations, see for instance Hirsch and Smale (1974). We get
$$x(t)=c_2{b - a + \la +\sqrt{\Delta}\over 2a} \exp(\nu_1 t)-c_1{a - b - \la +\sqrt{\Delta} \over 2a}\exp(\nu_2 t)$$
$$y(t)=c_1\exp(\nu_2 t) + c_2\exp(\nu_1 t).$$
We now find the the constants $c_1$ and $c_2$ by using the initial conditions $x(0)=0$ and $y(0)=1$. We have
$$c_1={b-a+\la+\sqrt{\Delta}\over 2\sqrt{\Delta}}\hbox{ and }c_2=1-c_1.$$
Since $\sqrt{\Delta}>|a+b -\la|$ we have $c_1>0$. It is also easy to check that $c_1<1$. Hence, $c_2>0$. 

It is not difficult to see that the function $y$ drops from $y(0)=1$ to some minimum and then increases to infinity as $t$ goes to infinity. Hence, there is a unique $T_c>0$ such that $y(T_c)=1$. The critical value $T_c$ is computed numerically by solving the equation $y(t)=1$. For any $T>T_c$ the Galton-Watson process $Z_k$ is super-critical and survives forever with a positive probability. For $T\leq T_c$ the process $Z_k$ dies out.  This completes the proof of Theorem 1.
\bigbreak

{\bf 4. Proof of Theorem 2.}

This proof is similar to the proof of Theorem 1. We now define a process $Z'_k$ that will play a role analogous to the one played by $Z_k$. 

We wait until the first killing time $T_1$ and we let $Z'_0$ be the number of bacteria in state 2. If $Z'_0=0$ we set $Z'_i=0$ for $i\geq 1$. If $Z'_0\geq 1$ then we wait until the second killing time $T_2$ and let $Z'_1$ be the number of  state 2 bacteria. More generally, for any $k\geq 1$ let $Z'_k$ be the number of state 2 individuals at the $(k+1)$-th  killing time $T_{k+1}$. 

The process $(Z_k')_{k\geq 1}$ is a branching process in an i.i.d. environment. To see this consider the 2's present at time $T_1+\dots+T_{k+1}$ for a fixed $k\geq 0$. Each one of these 2's is starting a patch that will last for the random time $T_{k+2}$. It is this random time that determines the random environment for each patch. Given $T_{k+2}$ all patches are independent and identically distributed. Moreover, the patch distribution does not depend on $k$. This fits the definition of a branching process in an i.i.d. environment. Smith and Wilkinson (1969) introduced branching processes in i.i.d. environments, see also Athreya and Karlin (1971). 
We will thus use the known criteria for survival and extinction for a branching process in an i.i.d. environment to prove Theorem 2. Note that so far everything we wrote is true for any distribution for the killing times. We start now concentrating on exponential killing times. 

Let $Y(T_1)$ be the number of 2's present at the first killing time $T_1$ starting the model with a single 2 at time 0.  Note that given $T_1=t$ the mean value of $Y(T_1)$ is $y(t)$
where $y(t)$ is the expected number of 2's given that the killing time $T_1>t$. The type of killing time (deterministic or random) does not affect the function $y$ and  we have already 
shown that
$$y(t)=c_1\exp(\nu_2 t) + c_2\exp(\nu_1 t).$$
Let us recall the criteria for the survival of the branching process in a random environment, in our context:

\bigskip
\noindent
{\bf Theorem (Smith and Wilkinson, 1969). } Assume that $E|\ln(y(T_1))|<\infty$.

a) If $E(\ln(y(T_1))) \le 0$, then extinction occurs with probability $1$.

b) If $E(\ln(y(T_1))) > 0$ and $E|\ln P(Y(T_1) \ge 1)| < \infty$, then the process survives with positive probability.

\bigskip
$\bullet$
We first check the integrability condition $E|\ln(y(T_1))|<\infty$.
In our setting, 
$$E|\ln(y(T_1))|=\int_0^\infty |\ln y(t)|\delta\exp(-\delta t)dt.$$
Recall that there is a positive real $T_c$ such that $y(t)<1$ for $t<T_c$ and $y(t)>1$ for $t>T_c$. We have
$$\int_{T_c}^\infty |\ln y(t)|\delta\exp(-\delta t)dt=\int_{T_c}^\infty \ln y(t)\delta\exp(-\delta t)dt.$$
Note that for $t\geq 0$, 
$y(t)\leq (c_1+c_2)\exp(\nu_1 t)=\exp(\nu_1 t)$.
Hence, 
$$\int_{T_c}^\infty |\ln y(t)|\delta\exp(-\delta t)dt\leq {\nu_1\over\delta}.$$
Hence, $E|\ln(y(T_1))|<\infty$ for $\delta>0$.

The next task is to study the sign of $m'$, where
$$m'=E(\ln(y(T_1)))=\int_0^\infty \ln (y(t))\delta\exp(-\delta t)dt.$$

\medskip
$\bullet$ We first show that $m'$ 
is positive for $\delta$ near 0.
Observe that 
$y(t)>c_2\exp(\nu_1 t)$
and therefore
$$\int_0^\infty \ln y(t) \delta\exp(-\delta t)dt\geq 
\ln c_2+{\nu_1\over\delta}.$$
Hence, for $\delta<-\nu_1/\ln c_2$ we have $m'>0$.
This ensures survival provided the integrability condition $E|\ln(P(Y(T_1)\geq 1)|<\infty$ holds. We now check this.
Recall that we start the process with a 2. If by time $T_1$ the 2 has not flipped to state 1 we have $Y(T_1)=1$. Hence,
$$P(Y(t)\geq 1)\geq e^{-bt} \quad \text{ and } \quad |\ln P(Y(t)\geq 1)|\leq bt.$$
Therefore,
$$E|\ln(P(Y(T_1)\geq 1)|=\int_0^\infty |\ln P(Y(t)\geq 1)|\delta\exp(-\delta t)dt\leq
\int_0^\infty bt\delta\exp(-\delta t)dt={b\over\delta}<\infty.$$
This proves the process survives with positive probability for $\delta$ small enough.

\bigbreak
$\bullet$ We now turn to the proof that $m'$ is negative
for $\delta$ large enough. By a change of variable we get
$$m'=\int_0^\infty \ln y({u\over\delta})\exp(-u)du.$$
The function $y$ is infinitely differentiable everywhere and a Taylor expansion gives
$$ y({u\over\delta})=y(0)+{u\over\delta}y'(0)+({u\over\delta})^2y''(c(u,\delta))$$
where $0<c(u,\delta)<{u\over\delta}$. Recall that $y(0)=1$  and that $\ln (1+x)\leq x$ for $x>-1$.
Hence,
$$\ln y({u\over\delta})\leq  {u\over\delta}y'(0)+({u\over\delta})^2y''(c(u,\delta)),$$
and
$$m'\leq y'(0)\int_0^\infty {u\over\delta}e^{-u}du+\int_0^\infty ({u\over\delta})^2y''(c(u,\delta))du=
{y'(0)\over\delta}+{1\over\delta^2}\int_0^\infty u^2y''(c(u,\delta))du.$$
Since $y'(0)<0$ to show that $m'<0$ for large $\delta$ it is enough to show that $\int_0^\infty u^2y''(c(u,\delta))du$ is bounded for $\delta$ larger than some $\delta_2$. We do so now. Note that since $\nu_1>\nu_2$
$$y''(t)=c_1\nu_2^2\exp(\nu_2 t) + c_2\nu_1^2\exp(\nu_1 t) \leq K\exp(\nu_1 t),$$
where $K=2\max(c_1\nu_2^2,c_2\nu_1^2)$.
Since $\nu_1>0$,
$$y''(c(u,\delta))\leq K\exp(\nu_1c(u,\delta) )\leq K\exp(\nu_1{u\over\delta}).$$
Hence, for $\delta>\nu_1$
$$0\leq \int_0^\infty u^2y''(c(u,\delta))du\leq K \int_0^\infty u^2\exp(\nu_1{u\over\delta})du=
{2K\over (1-\nu_1/\delta)^3}.$$
Therefore,
$$0\leq \limsup_{\delta\to\infty}\int_0^\infty u^2y''(c(u,\delta))du\leq 2K.$$
This proves extinction for $\delta$ large enough.

\medskip
$\bullet$ We now prove that the probability of surviving is a non-increasing function of the intensity $\delta$ of killing times.

In the specific case of exponential killing times, there exists a natural coupling of the processes for distinct values of the killing parameter $\delta$ that we now explain. Assume that the delays $(T_k)_{k \ge 1}$ between two killing times are independent exponential random variables with parameter $\delta$, then the sequence of killing times $(S_k=\sum_{i=1}^k T_i)_{k \ge 1}$ is a Poisson point process with intensity $\delta$. The classical decimation procedure allows to make a coupling between two Poisson point processes $(S_k)_{k \ge 1}$ and $(S'_k)_{k \ge 1}$ with respective intensities $\delta'>\delta$ such that $(S_k)_{k \ge 1}$ is a subsequence of $(S'_k)_{k \ge 1}$. 

Let us explain briefly this procedure. Let $(T'_k)_{k \ge 1}$ be an i.i.d sequence of exponential random variables with parameter $\delta'$: the killing times $(S'_k=\sum_{i=1}^k T'_i)_{k \ge 1}$ are a Poisson point process with parameter $\delta'$. Consider next a sequence $(\varepsilon_k)_{k \ge 1}$ of i.i.d. Bernoulli random variables with parameter $\delta/\delta'$. Then in the sequence $(S'_k)_{k \ge 1}$, we keep the $k$-th point $S'_k$ if $\varepsilon_k=1$ and erase it otherwise. Doing so, we obtain a new sequence $(S_k)_{k \ge 1}$  that is a Poisson point process with intensity $\delta$. 

Let us now explain a graphical construction for our bacteria growth process.
\begin{itemize}

\item Let $\mathbb T$ be the following binary tree. The set $V$ of vertices is the set of finite words on the alphabet $\{0,1\}$, the root is the empty word $\emptyset$, and a non-empty word $v_1v_2\dots v_{n-1}v_n$ is linked to his parent $v_1v_2\dots v_{n-1}$ and to his two children $v_1v_2\dots v_{n-1}v_n0$ and $v_1v_2\dots v_{n-1}v_n1$. 

\item We now add a time coordinate. Consider a family $(S_v)_{v \in V}$ of iid random variables following the exponential law with parameter $\lambda$: $S_v$ represents the time between the birth of the particle at $v$ and its splitting time (i.e. when the particle gives birth). Set
$$T_\emptyset =0 \text{ and }\forall v \in V\backslash \{\emptyset \} \quad T_v= \sum_w S_w,$$
where the sum is taken on all the ancestors $w$ of $v$ ($\emptyset$ included, $v$ excluded).  

We  build a new tree $\tilde{\mathbb T}$ in $\mathbb T \times [0,+\infty)$: the set $\tilde V$ contains the points
$$ (v, T_v)\text{ and } (v, T_v+S_v) \text{ for } v \in V.$$
For each $v \in V$ we do the following.
We draw a vertical edge between $(v, T_v)$ and $(v, T_v+S_v)$: this edge corresponds to the time between the birth at $v$ and its splitting time. We also draw an horizontal edge between $(v, T_v+S_v)$ and each of its two children $(v\epsilon, T_v+S_v)$, $\epsilon \in \{0,1\}$: this edge corresponds to the splitting of the particle $v$ in two. 
Let $\eta_t$ be the number of living particles at time $t$, or, equivalently in the graphical setting, the number of vertical edges intersecting the horizontal plane  with time coordinate $t$. Then $(\eta_t)_{t \ge 0}$ is a continuous time branching process with rate $\lambda$.

\item Now we add the switches between the two states $1$ and $2$. Consider, independently of everything built before, two families of independent Poisson point processes $(\omega^1_v)_{v \in V}$ and $(\omega^2_v)_{v \in V}$ such that $\omega^1_v$ and $\omega^2_v$ are independent Poisson point processes on the half-line $\{v\}\times [0,+\infty)$ with respective intensity $a$ and $b$. Now we color the tree in white and red as follows. We start from $(\emptyset,0)$ in white, we color each branch of the tree following the time direction, changing the color to red every time we meet a point in one on the Poisson point processes $(\omega^2_v)_{v \in V}$, changing the color to white every time we meet a point in one on the Poisson point processes $(\omega^1_v)_{v \in V}$. White segments correspond with bacteria in normal state, while red segments correspond to persistent bacteria.

\item As persistent bacteria cannot split, at each red splitting point, we cut the horizontal red edge that goes to the son whose name ends with a $1$ and erase the corresponding subtree. 

\end{itemize}

Then the number
$\eta^1_t$ of white vertical edges -- resp. $\eta^2_t$ of red vertical edges -- intersecting the horizontal plane  with time coordinate $t$ correponds to the number of normal -- resp. persistent -- bacteria living at time $t$ in our model when no killing occurs.

Independently, we consider the previous coupling of the two Poisson point processes $(S_k)_{k \ge 1}$ and $(S'_k)_{k \ge 1}$ of killing times. To recover the bacteria process with killing rate $\delta$ -- resp. $\delta'$ --, we just have to cut, at every time $S_k$ -- resp $S'_k$ -- all the white edges, and to erase the corresponding subtrees. As  $(S_k)_{k \ge 1}$ is a subsequence of $(S'_k)_{k \ge 1}$, it is clear that the bacteria process with killing intensity $\delta'$ is a subtree of the bacteria process with killing intensity $\delta$. In particular, the probability of surviving is a non-increasing function of $\delta$.

\medskip
$\bullet$ We finally remark that $\delta \mapsto m'(\delta)=\int_0^\infty \ln (y(t))\delta\exp(-\delta t)dt$ is a continuous function of $\delta$, which ensures that $m'(\delta_c) \le 0$, and thus that the process almost surely dies out when $\delta=\delta_c$.

\bigbreak

{\bf References}

K.B. Athreya and S. Karlin (1971). On branching processes with random environments: I extinction probabilities. The Annals of Mathematical Statistics 42, 1499-1520.

N.Q. Balaban, J. Merrin,R. Chait and S. Leibler (2004). Bacterial persistence as a phenotypic switch. Science 305, 1622-1625.

J.W. Bigger (1944). Treatment of staphylococcal infections with penicillin. Lancet, 497-500.

D. Dubnau and R. Losick (2006). Bistability in bacteria. Mol. Microbiol. 61, 564-572.

A. Gardner, S. A. West and A.S. Griffin (2007). Is bacterial persistence a social trait? PLos One 2(8): e752.

P.L. Graumann (2006). Different genetic programmes within identical bacteria under identical conditions: the phenomenon of bistability greatly modifies our view on bacterial populations. Mol. Microbiol. 61, 560-563.

M. W. Hirsch and S. Smale (1974) {\sl Differential equations, dynamical systems and linear algebra}. Academic Press.

E. Kussell, R. Kishony, N.Q. Balaban and S. Leibler (2005). Bacterial persistence: a model of survival in changing environments. Genetics, 169, 1807-1814.

B.R. Levin (2004). Noninherited resistance to antibiotics. Science 305, 1578-1579.

B.R. Levin and Rozen  (2006). Non-inherited antibiotic resistance. Nature Rev. Microbiol. 4, 556-562.

K. Lewis (2007). Persister cells, dormancy and infectious disease. Nat. Rev. Microbiol. 5, 48-56.

R. B. Schinazi (1999). {\sl Classical and spatial stochastic processes.} Birkhauser.

W.L. Smith and W.E. Wilkinson (1969). On branching processes in random environments. The Annals of Mathematical Statistics 40, 814-827.

C. Wiuff, R.M. Zappala, R.R. Regoes, K.N. Garner, F. Baquero and B.R. Levin (2005). Phenotypic tolerance: antibiotic enrichment of noninherited resistance in bacterial populations. Antimicrobial agents and chemotherapy 49, 1483-1494.

\end{document}